# Adaptive Non-singular Terminal Sliding Mode Fault-tolerant Control of Robotic Manipulators Based on Contour Error Compensation

**Zhu Dachang**[1,*], **Du Baolin**[1], **Cui Aodong**[1], **Zhu Puchen**[2]

1. School of Mechanical and Electrical Engineering, Guangzhou University, Guangzhou 510006, China
2. School of Automation, Guangdong University of Technology, Guangzhou 510006, China

* Correspondence: zdc98998@gzhu.edu.cn; Tel.: +86-020-3936-6923

**Abstract:** To achieve accurate contour tracking of robotic manipulators with dynamic uncertainties, coupling and actuator faults, an adaptive non-singular terminal sliding mode control (ANTSMC) based on cross-coupling is proposed. Firstly, the singularity is eliminated completely by using a terminal sliding mode manifold. Secondly, an adaptive tuning approach is selected for avoid the demand of the bound of system uncertainty, and the stability of the proposed control strategy is demonstrated by the sense of the finite-time stability theory. Furthermore, the cross-coupled ANTSMC law is proposed for contour tracking at the end-effectors level of robotic manipulators. Thirdly, a unified framework of cross-coupling contour compensation and reference position pre-compensation is designed by combining cross-coupling control with parabolic transition trajectory planning. Finally, numerical simulation and experimental results are shown to prove the effectiveness of the proposed control strategy.

## 1. Introduction

During recent year, with the rapid development of modern industrial technology, robot have been diffusely used in handling, assembling, welding and other fields. In these applications, trajectory tracking with high precision and fast dynamical response is necessarily required in order to accomplish specific operations [1]. However, due to the long-term operation of robotic manipulators and the wear and aging of the joint parts, robotic manipulators are prone to various unpredictable disturbances, which will affect the trajectory tracking accuracy and reliability of robotic manipulators directly. It is challenging to fulfill accurate and fast trajectory tracking and reliability for robotic manipulators due to the impediments in designing valid trajectory tracking controllers. The primary impediments involve the influence of dynamic uncertainties, strongly coupled and actuator faults and so on [2-3].

Sliding mode control (SMC) is a particular and powerful class of variable structure control essentially [6], and the structure of control system is also not fixed, but change according to the current state of the system in the dynamic process, which forces the system state to follow the predetermined sliding mode manifold trajectory, and the sliding mode manifold can be designed and it is independent of target parameters and external disturbances [8]. Due to the above superiorities, SMC has been diffusely used in the precise trajectory tracking control of uncertain robotic manipulators [10-13]. However, the fast convergence of system states to the equilibrium point is one of the important performance indicators in SMC. But these traditional SMC can guarantee only asymptotical convergence of system states because of the asymptotical convergence of the linear sliding mode manifold that are generally selected [16]. To solve this problem, the TSMC with global finite-time stability was proposed by inducting nonlinear function in the traditional SMC [17]. The finite-time control can provide fast transient and high precision performance [18]. However it does not resolve the singularity phenomenon near the equilibrium point caused by the negative exponent of the TSMC. In order to eliminate the singularity completely for the global finite-time tracking problem of robotic manipulator with dynamic uncertainties and external disturbance, an integral sliding mode manifold and its TSMC was proposed by Su [14], and manifested the global finite-time convergence of both the sliding mode manifold and tracking error.

Although SMC offers a valid solution to accurate trajectory tracking of uncertain robotic manipulator because of the strong robustness with respect to uncertainty [14]. However the traditional SMC relies on a priori information of the frontier of system uncertainties in the design process, and this assumption is difficult to evaluate the frontier information of uncertainties such as random fault parameters and disturbances in advance practically. In order to solve this problem, the switching gain is designed too high to reduce the influence of uncertainty estimation on the control system usually [20, 21], but it also causes the phenomenon of high frequency chattering. To solve the problem of the unknown upper bound of system uncertainty in practical tasks, a new ANTSMC was proposed by Boukattaya [22]. In further work, a fraction-order adaptive back-stepping approach for robot in the existence of actuator faults was investigated by Anjum [23], and the controller adopted a new adaptive technology because of which the controller does not need the fault and uncertainty information existing in the system.

At present, the trajectory tracking control of robotic manipulator focus on optimizing the tracking performance of each joint, so as to improve the contour machining accuracy of the end-effector indirectly. However, the problem of contour error degradation caused by lack of coordination of each joint axis of robotic manipulator is ignored [24]. Moreover, the joint axis of robotic manipulator with fast response cannot be fully utilized to compensate the influence of the joint axis with large position tracking error on the overall contour control performance in real time. Therefore, the above approach of only improving the trajectory tracking accuracy of each joint alone cannot resolve the problem of profile machining accuracy of the end-effector of robotic manipulator with multi-axis linkage fully. For the high precision contour control of CNC machine tools with orthogonal structure, a cross-coupling control method was proposed by Reason [25], and the coordination of each axis of the machine tool was improved by compensating each axis with contour error, so as to improve the contour accuracy of the workpiece. For studying the contour tracking control of serial robotic manipulator, a cross-coupling PD controller based on contour error of end-effector and tracking error of joint was proposed by Ouyang [26], so as to ensure the contour tracking accuracy of the end-effector of robotic manipulator. However, there are few investigations on cross-coupling contour tracking control for the factors of the non-smooth path. Due to the fact that contour motion control is related to the motion trajectory closely, the sudden change of position, velocity and acceleration at the inflection point of the non-smooth path may lead to the resonance problem of the mechanism.

For details, in order to further optimize the contour tracking property of robot, a new cross-coupled ANTSMC law is proposed, and it is a strategy combining the ANTSMC law for trajectory tracking at joint level and the PD control law for contour tracking at the end-effectors level of robotic manipulators. While the proposed ANTSMC can improve the trajectory tracking accuracy of each joint of the robotic manipulators with dynamic uncertainties and actuator faults, coupling factors are introduced among the multi-axes to enhance the coordination control property of system by cross-coupled control. An unified framework of cross-coupling contour compensation and reference position pre-compensation is designed by the cross-coupling technology combined with Cartesian space trajectory planning technology to introduce the necessary coordination mechanism and then from the wider research and realization of fast and high accurate movement of robot. The primary contributions of this paper are summarized as followings.

(1) To eliminate the singularity completely, a terminal sliding manifold is proposed. Moreover, adaptive tuning approach is adopted to compensate dynamic uncertainties and actuator faults during the operation of robotic manipulators, and the stability of the proposed control strategy is demonstrated by the sense of the finite-time stability theory.

(2) Combining the joint trajectory tracking with the end-effector contour tracking, a novel cross-coupled ANTSMC law is proposed.

(3) A unified framework of cross-coupling contour compensation and reference position pre-compensation is designed.

This paper is organized as follows: Section 2 indicates the problem formulation and motivation. Section 3 presents the design of ANTSMC for the precise trajectory tracking of robotic manipulators with dynamic uncertainties and actuator faults, and its finite time stability is discussed. Section 4 describes the contour error compensation of robotic manipulators based on cross coupling control. Numerical simulation and experiments results and conductions are presented in Section 5 and Section 6, respectively.

**2. Problem formulation and motivation**

The dynamic of $n-$DoF (Degree-of-Freedom) robotic manipulators can be expressed by Newton-Euler formula as [1]

$$M(q)\ddot{q}+B(q,\dot{q})\dot{q}+G(q)+F(q,\dot{q})=\tau-\tau_d \tag{1}$$

where $q,\dot{q},\ddot{q}\in R^{n\times 1}$ denotes the joint position, velocity, and acceleration vectors, respectively, $M(q)\in R^{n\times n}$ expresses the positive definite inertia matrix, $B(q,\dot{q})\in R^{n\times n}$ denotes the Coriolis, centripetal matrix, $G(q)\in R^{n\times 1}$ is the gravity vector, $\tau_d\in R^{n\times 1}$ represents external disturbance torque vector, $\tau\in R^{n\times 1}$ denotes the torque input vector, and $F(q,\dot{q})\in R^{n\times 1}$ denotes the friction terms.

For actual applications, it is difficult to obtain the precise dynamic model of the robotic manipulator because of the friction and external disturbance. Hence, the dynamic model (1) can be written as

$$M(q)\ddot{q}+B(q,\dot{q})\dot{q}+G(q)+\varphi(q,\dot{q},t)=\tau \tag{2}$$

where $\varphi(q,\dot{q})$ denotes the lumped uncertainty of the system and can be defined as

$$\varphi(q,\dot{q},t)=\Delta M(q)\ddot{q}+\Delta B(q,\dot{q})\dot{q}+\Delta G(q)+\tau_d+F(q,\dot{q}) \tag{3}$$

The robotic manipulator dynamics (2) satisfies the following properties:

**Property 1:** [1] $M(q)$ is the symmetric and positive definite matrix, and bounded:

$$mI<M(q)=M^T(q)\leqslant\bar{m}I \tag{4}$$

where $m,\bar{m}$ are positive constants, respectively, and $0<m<\bar{m}$.

**Property 2:** [2] The matrix $M(q)-2B(q,\dot{q})$ is skew symmetric and satisfies

$$D^T\left[M(q)-2B(q,\dot{q})\right]D=0 \tag{5}$$

where $D$ is any nonsingular matrix.

**Property 3:** [2] $G(q)$ is bounded by

$$\|G(q)\|\leqslant G_k \tag{6}$$

where $G_k$ is positive constant.

Further, consider the problem of actuator faults during the operating process of robotic manipulators, and the robot dynamics (2) can be written by [10]

$$M(q)\ddot{q}+B(q,\dot{q})\dot{q}+G(q)+\varphi(q,\dot{q},t)=\tau+f \tag{7}$$

where $f=\gamma(t-T_f)\phi(q,\dot{q},\tau)$ denotes the function of actuator faults, $\gamma(t-T_f)\in R^{n\times 1}$ is the time profile of the faults, $\phi(q,\dot{q},\tau)$ is the fault vector, and $T_f$ is the time of appearance of the faults.

The time profile of the faults $\gamma(t-T_f)$ is a diagonal matrix, and can be expressed as

$$\gamma(t-T_f)=diag\left[\gamma_1(t-T_f),\gamma_2(t-T_f),\cdots,\gamma_n(t-T_f)\right] \tag{8}$$

where $\gamma_i$ represents the influence of fault to the $i^{th}$ state equation.

The time profile mode of the fault can be given by

$$\gamma(t-T_f)=\begin{cases}0 \quad if \quad t<T_f \\ 1-e^{-\sigma_i}(t-T_f) \quad if \quad t\geq T_f\end{cases} \tag{9}$$

where $\sigma_i>0$ is the evolution rate of the fault.

The objective is that design an control strategy such that the trajectory tracking error of each joint of the robotic manipulators converge zero and improve the coordination of each joint, so that the end-effector can track the non-smooth path with the sudden change of position, velocity and acceleration well in the existence of both fault and uncertainty.

## 3. Adaptive non-singular terminal sliding mode control

3.1 Non-singular terminal sliding mode control (NTSMC) design

The position tracking error denoted by $e(t)\in R^{n\times 1}$ in joint space is defined as

$$e=q-q_d \tag{10}$$

where $q_d\in R^{n\times 1}$ is the expected trajectory.

A non-singular terminal sliding mode manifold $s$ is defined as

$$s=\dot{e}+c_1e^{\alpha/\beta}+c_2e^{\eta}\,sgn(e) \tag{11}$$

where $c_1,c_2\in R^{n\times n}$ are constant positive-definite diagonal matrices, respectively, $\alpha$ and $\beta$ are positive odd integers, respectively, which satisfy the condition: $1<\alpha/\beta<2$, $\eta>0$, $sgn(*)$ is the signum function.

The derivative of sliding mode manifold (11) can be derived as

$$\dot{s}=\ddot{e}+\frac{c_1\alpha}{\beta}e^{\alpha/\beta-1}\dot{e}+c_2\eta e^{\eta-1}\dot{e} \tag{12}$$

Considering the positive-definite Lyapunov function

$$V_1 = \frac{1}{2} s^T M s \tag{13}$$

The time derivative of $V_1$ and using (12) yields

$$\begin{aligned} \dot{V}_1 &= \frac{1}{2}\dot{s}^T M s + \frac{1}{2} s^T \dot{M} s + \frac{1}{2} s^T M \dot{s} = s^T M \dot{s} + \frac{1}{2} s^T \dot{M} s \\ &= s^T M\left(\ddot{q} - \ddot{q}_d\right) + s^T B s + s^T M\left(\frac{c_1 \alpha}{\beta} e^{\left(\alpha/\beta\right)-1}\dot{e} + c_2 \eta e^{\eta-1}\dot{e}\right) \end{aligned} \tag{14}$$

Substituting (7) into (14), one obtains

$$\dot{V}_1 = s^T M\left[M^{-1}\left(\tau - B\dot{q} - G - \varphi + f\right)\right] - s^T M \ddot{q}_d + s^T M\left(\frac{c_1 \alpha}{\beta} e^{\alpha/\beta - 1}\dot{e} + c_2 \eta e^{\eta-1}\dot{e}\right) + s^T B s \tag{15}$$

Simplifying (15) yields

$$\dot{V}_1 = s^T\left[\tau + B\left(s - \dot{q}\right) - G - \varphi + f + M\left(\frac{c_1 \alpha}{\beta} e^{\alpha/\beta - 1}\dot{e} + c_2 \eta e^{\eta-1}\dot{e} - \ddot{q}_d\right)\right] \tag{16}$$

By setting $\dot{V}_1 = 0$, the equivalent control law $u_{eq}$ is derived as

$$u_{eq} = -B\left(s - \dot{q}\right) + G + \varphi - f - M\left[\frac{c_1 \alpha}{\beta} e^{\alpha/\beta - 1}\dot{e} + c_2 \eta e^{\eta-1}\dot{e} - \ddot{q}_d\right] \tag{17}$$

The equivalent control law $u_{eq}$ can force the system state remain on the sliding mode manifold if the dynamic of system is known accurately. However, this assumption is difficult to satisfy in factual applications. For appease the sliding condition in the presence of uncertainties, a switching control $u_{sw}$ is introduced as follows:

Assumed that the actuator faults and the lumped uncertainty of the system are denoted by function $K$, yields to

$$K = \varphi - f \tag{18}$$

Then, the upper bound of function $K$ is estimated as follows

$$\left\|K\right\|_{max} = \left\|\varphi\right\| + \left\|f\right\| \tag{19}$$

where $\left\|*\right\|$ is the standard Euclidean norm.

The switching control law $u_{sw}$ is given by

$$u_{sw} = -ks\left|s\right| - k\,sgn\left(s\right)\left|s\right| \tag{20}$$

where $k = \left\|K\right\|_{max} + \upsilon$, $\upsilon$ is the switching control gain, and $\upsilon \geq 0$.

Thus, the NTSMC is derived as

$$
\begin{aligned}
u &= u_{eq} + u_{sw} \\
&= -M\left[\frac{c_1\alpha}{\beta} e^{\alpha/\beta - 1}\dot{e} + c_2\eta e^{\eta-1}\dot{e} - \ddot{q}_d\right] - B(s - \dot{q}) + G + K - \left(\|K\|_{max} + \upsilon\right)\left(s|s| + sgn(s)|s|\right)
\end{aligned}
\tag{21}
$$

3.2 Adaptive non-singular terminal sliding mode control (ANTSMC) design

The design process of NTSMC relies on the upper bound of function $K$. However, it is difficult to appraise the upper bound of random fault parameters and the lumped uncertainty of the system in practical application. In this section, an ANTSMC is designed to appraise the unknown upper bounds of the presence of both fault and uncertainty.

The estimation error is defined as $\tilde{K} = K - \hat{K}$, $\hat{K}$ is the estimated value of $K$. Suppose the uncertainty of robotic manipulator system change slowly, that is $\dot{K} = 0$. The time derivative of the estimation error $\tilde{K}$ can be derived as follows

$$
\dot{\tilde{K}} = \dot{K} - \dot{\hat{K}} = -\dot{\hat{K}} \tag{22}
$$

Considering the positive-definite Lyapunov function

$$
V_2 = V_1 + \frac{1}{2\xi}\tilde{K}^T\tilde{K} \tag{23}
$$

where $\xi$ is positive integer.

The time derivative of (23) is obtained as follows

$$
\begin{aligned}
\dot{V}_2 &= \dot{V}_1 + \frac{1}{\xi}\tilde{K}^T\dot{\tilde{K}} \\
&= s^T\left[\tau + f - B\dot{q} - G - \varphi - M\ddot{q}_d + Bs + M\left[\frac{c_1\alpha}{\beta} e^{\alpha/\beta - 1}\dot{e} + c_2\eta e^{\eta-1}\dot{e} - \ddot{q}_d\right]\right] - \frac{1}{\xi}\tilde{K}^T\left(\dot{\hat{K}} - \xi\left(M^{-1}\right)^T s\right)
\end{aligned}
\tag{24}
$$

Adaptive law is given by

$$
\dot{\hat{K}} = \xi\left(M^{-1}\right)^T s \tag{25}
$$

Then, ANTSMC can be designed as follows

$$
u = -M\left[\frac{c_1\alpha}{\beta} e^{\alpha/\beta - 1}\dot{e} + c_2\eta e^{\eta-1}\dot{e} - \ddot{q}_d\right] - B(s - \dot{q}) + G + \hat{K} - \hat{k}s|s| - \hat{k}\,sgn(s)|s| \tag{26}
$$

The schematic for the proposed ANTSMC strategy is displayed in Fig.1.

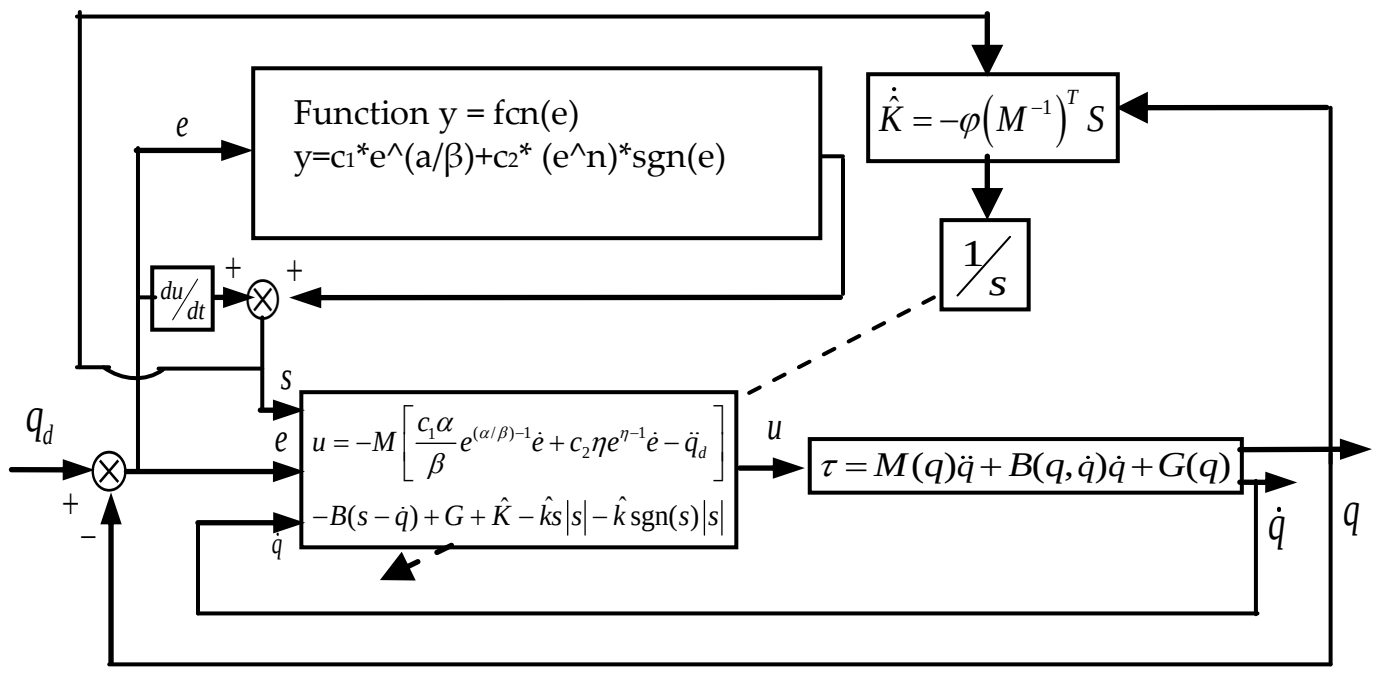

Fig.1 Schematic of the proposed ANTSMC strategy.

3.3 Stability Analysis

*Theorem* 1: For finite-time stability with fast time convergence, the Lyapunov function $V(x)$ with initial value $V_0$ is given as

$$\dot{V}(x)+aV(x)+bV^{\delta}(x)\le 0,\quad \forall x\ge x_0,\quad V(x_0)\ge 0 \tag{27}$$

where $a>0$, $b>0$, $0<\delta<1$, $V(x)$ satisfies the inequality at any $x_0$.

Thus, the matching stability time $T$ can be derived as

$$T\le\frac{1}{a(1-\delta)}In\frac{aV^{1-\delta}(x_0)+b}{b} \tag{28}$$

**Proof:** Consider the following Lyapunov function:

$$V=\frac{1}{2}s^{T}Ms \tag{29}$$

The time derivative of $V$ can be derived as

$$\dot{V}=s^{T}M\dot{s}+\frac{1}{2}s^{T}\dot{M}s \tag{30}$$

Substituting (12) into (30), one obtains

$$\begin{aligned}\dot{V}&=s^{T}M\left(\ddot{e}+\frac{c_1\alpha}{\beta}e^{(\alpha/\beta)-1}\dot{e}+c_2\eta e^{\eta-1}\dot{e}\right)+s^{T}Bs\\&=s^{T}\left[\tau+f-B\dot{q}-G-\varphi+Bs+M\left(\ddot{e}+\frac{c_1\alpha}{\beta}e^{(\alpha/\beta)-1}\dot{e}+c_2\eta e^{\eta-1}\dot{e}-\ddot{q}_d\right)\right]\end{aligned} \tag{31}$$

By putting (26) in (31) yields

$$\dot{V}=s^{T}\left[-\tilde{K}-\hat{k}s|s|-\hat{k}\,sgn(s)|s|\right] \tag{32}$$

Simplifying (32), and implying $s^{T}s=\|s\|^{2}$ and $s^{T}\,sgn(s)=\|s\|$, the above equation (32) can be derived as

$$\dot{V}=-\|\tilde{K}\|\|s\|-\hat{k}\|s\|^{2}\|s\|-\hat{k}\|s\|\|s\| \tag{33}$$

Assume the function of actuator faults and the lumped uncertainty of the system $K$ changes slowly, the above equation (33) can be derived as

$$\dot{V} = -\left(\hat{k}\|s\|^2 + \hat{k}\|s\|\right)\|s\| \le -\left(s^T\hat{k}s + s^T\hat{k}\,sgn(s)^{\mu}\right)\|s\| \tag{34}$$

where $\mu \cong 1$.

As

$$s^T ks^r = \sum_{i=1}^{n} k_i s_i^{r+1} \ge \gamma\left\{\sum_{i=1}^{n}\frac{1}{2}\bar{m}s_i^2\right\}^{\lambda} \ge \gamma\left(\frac{1}{2}s^T Ms\right)^{\lambda} \tag{35}$$

where $\lambda \triangleq {(1+r)}/{2}$, $\gamma \triangleq k_{min}\left\{{2}/{\bar{m}}\right\}^{\lambda}$, $k_{min} \triangleq min(k_i)$.

Then

$$\dot{V} \le -\left(s^T\hat{k}s + s^T\hat{k}\,sgn(s)^{\mu}\right)\|s\| \le -\|s\|\gamma\left(\frac{1}{2}s^T Ms\right) - 1^{{(\mu+1)}/{2}}\|s\|\gamma\left(\frac{1}{2}s^T Ms\right)^{{(\mu+1)}/{2}} \tag{36}$$

Simplifying (36), one can obtain

$$\dot{V} \le -\|s\|\gamma V - 1^{{(\mu+1)}/{2}}\|s\|\gamma V^{{(\mu+1)}/{2}} \tag{37}$$

According to theorem 1, the proposed adaptive nonsingular terminal sliding mode controller is finite time stable.

When $x_0 = 0$, the stable time can be expressed as

$$t_s \le \frac{1}{\Omega_1(1-\delta)} In\left(1 + \frac{\Omega_1 V_0^{{(1-\mu)}/{2}}}{\Omega_2}\right) \tag{38}$$

where $\Omega_1 = \|s\|\lambda_{max}\left(\hat{k}\right)$, $\Omega_2 = 1^{{(\mu+1)}/{2}}\|s\|\lambda_{max}\left(\hat{k}\right)$.

## 4. Contour Error Compensation of Robot Based on Cross Coupling Control

Since the end-effector of the robotic manipulators track time-varying points rather than spatial contour, it may lead to the end-effector leaving the preset machining contour line to catch up with the reference position specified by the trajectory at the sampling time, resulting in phenomenon of radial contraction [15]. In this section, combining the joint trajectory tracking with the end-effector contour tracking, a novel cross-coupled ANTSMC is proposed.

### 4.1 Contour Error of End-Effector

The trajectory planning interpolation method can be used to fit the trajectory contour into a straight line or a circular contour. The contour error model of the straight line contour is displayed in Fig.2.

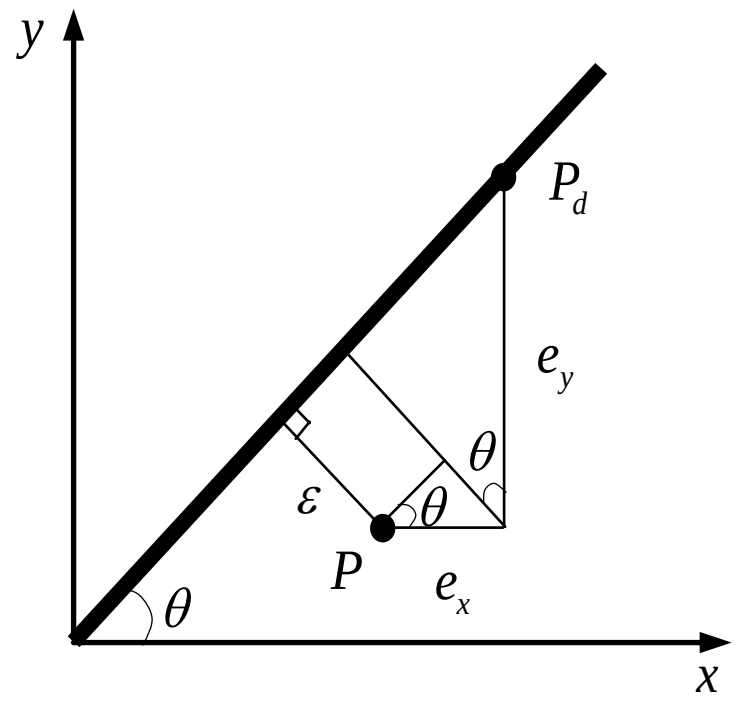

Fig.2 The contour error model of the straight line contour

In Fig.2, $P$ is the actual position of the end-effector of robotic manipulator, $P_d$ is the reference point of design, $\varepsilon$ is the contour error, $e_x$ and $e_y$ are the $x$ axis and $y$ axis error components of the tracking error $e$ of end-effector, respectively. $\theta$ denotes an angle between the reference trajectory and the $x$ axis.

Defined the contour error as the shortest distance between the prevailing position and the desired contour curve [26], the contour error of plane line can be defined as

$$\varepsilon = c_x e_x + c_y e_y \tag{40}$$

where $c_x = -sin\theta$, $c_y = cos\theta$ are the cross coupling operator.

## 4.2 Adaptive Non-singular Terminal Sliding Mode Control (ANTSMC) Based on Cross Coupling

The cross coupling method is used to design the contour motion controller as follows

$$u_i = K_p \varepsilon + K_d \dot{\varepsilon} \tag{41}$$

where $K_p$ and $K_d$ denote the proportional gain and the differential gain of the contour motion controller, respectively.

Combined with (20), the ANTSMC based on cross coupling can be derived as

$$u = -M\left[\frac{c_1\alpha}{\beta} e^{\left(\alpha/\beta\right)-1}\dot{e} + c_2\eta e^{\eta-1}\dot{e} - \ddot{q}_d\right] - B(s-\dot{q}) + G + \hat{K} - \hat{k}s|s| - \hat{k}\,sgn(s)|s| + \left(K_p\varepsilon + K_d\dot{\varepsilon}\right)C_n \tag{42}$$

where $C_n \in R^{n\times 1}$ is the contouring error compensation rectifier gain, which is used to correct the joint component of contouring error compensation.

The schematic of the proposed ANTSMC strategy based on cross coupling is shown in Fig.3.

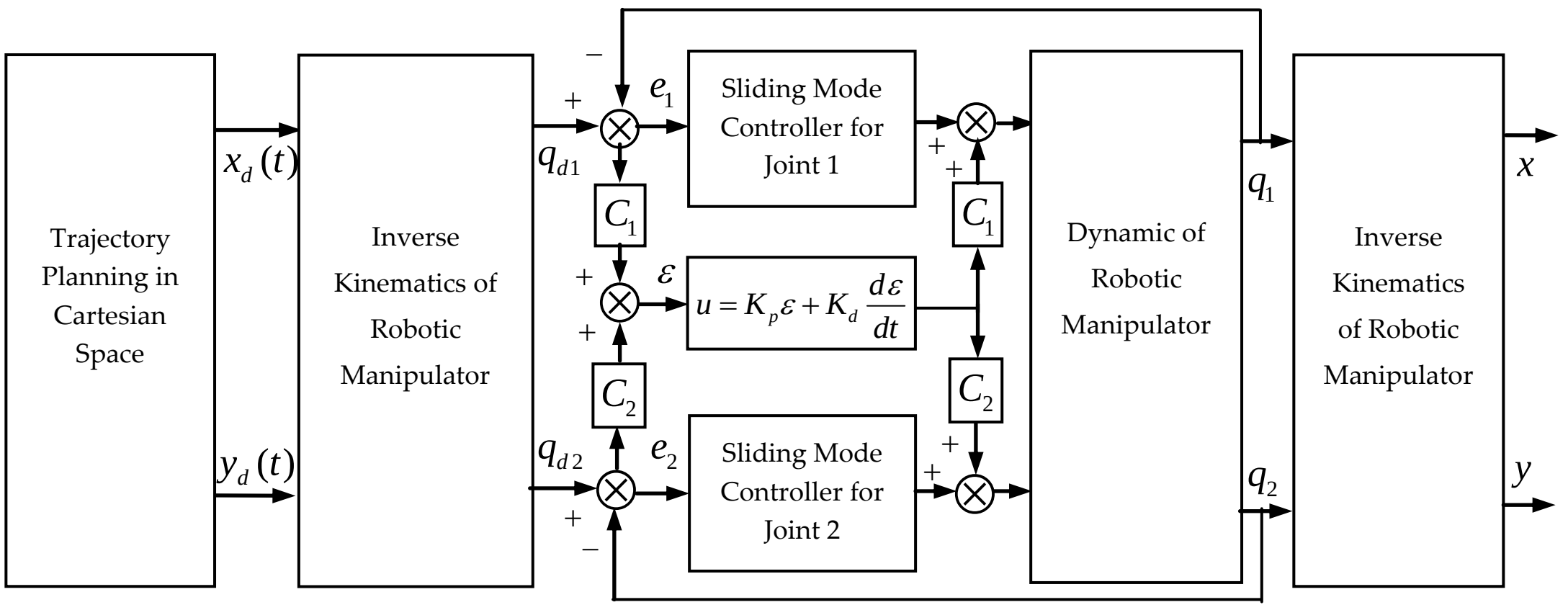

Fig.3 Schematic of the proposed ANTSMC strategy based on cross coupling.

## 5. Result and Discussions

In this section, numerical simulation and experiments are shown out for verify the effectiveness of the proposed control method. Simulations are divided into two parts, simulation of path tracking without planning, simulation of path tracking via parabola transition planning, are applied in the simulations with the intention of comparing the contour tracking accuracy and robustness against dynamic uncertainties and actuator faults. Experiments results are presented in another parts.

### 5.1. Simulation Results

To illustrate the effectiveness of the proposed control method, consider a Two-link robotic manipulator plant displayed in Fig 4.

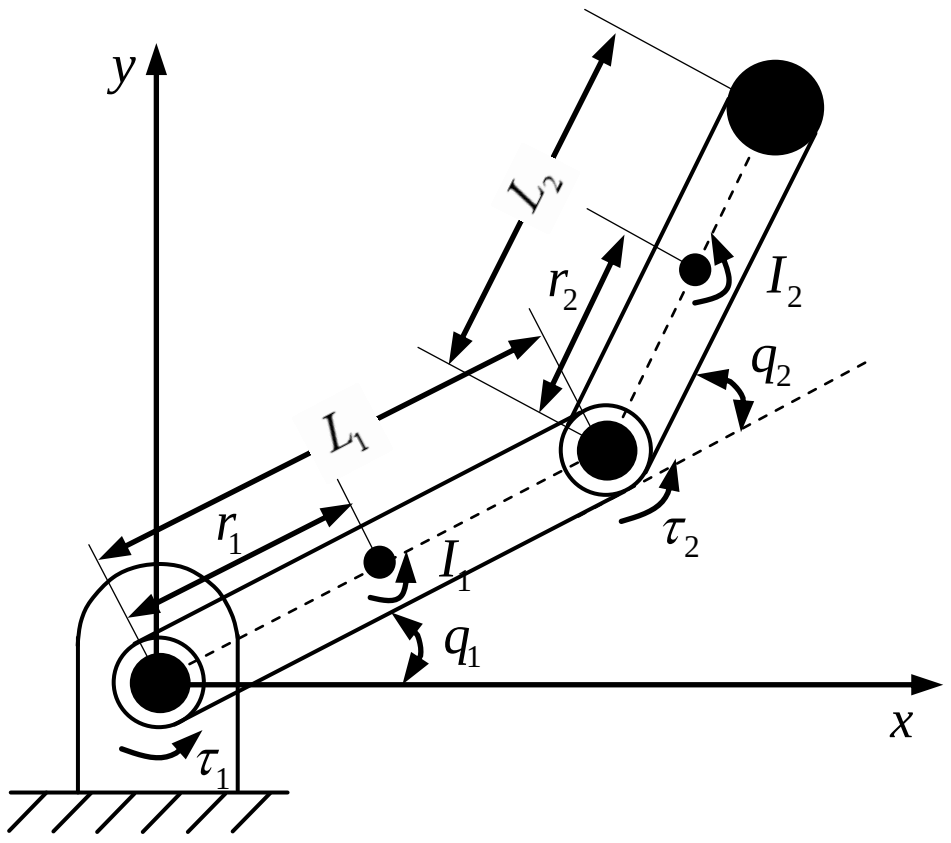

Fig.4 Two-link robotic manipulator plant.

Assumed that the mass of each link is concentrated. The dynamic of the robotic manipulator can be expressed by

$$\begin{bmatrix} \tau_1 \\ \tau_2 \end{bmatrix} = \begin{bmatrix} M_{11} & M_{12} \\ M_{21} & M_{22} \end{bmatrix} \begin{bmatrix} \ddot{q}_1 \\ \ddot{q}_2 \end{bmatrix} + \begin{bmatrix} B_{11} & B_{12} \\ B_{21} & B_{22} \end{bmatrix} \begin{bmatrix} \dot{q}_1 \\ \dot{q}_2 \end{bmatrix} + \begin{bmatrix} G_1 \\ G_2 \end{bmatrix} \tag{43}$$

where $[\tau_1 \quad \tau_2]^T$ denotes the torque input vector, $M_{11} = m_1 r_1^2 + m_2 r_2^2 + m_2 L_1^2 + 2m_2 L_1 r_2 \cos q_2 + I_1 + I_2$, $M_{12} = m_2 r_2^2 + m_2 L_1 r_2 \cos q_2 + I_2$, $M_{21} = m_2 r_2^2 + m_2 L_1 r_2 \cos q_2 + I_2$, $M_{22} = m_2 r_2^2 + I_2$, $B_{11} = -2m_2 L_1 r_2 \sin q_2 \dot{q}_2$, $B_{12} = -2m_2 L_1 r_2 \sin q_2 \dot{q}_2 - 2m_2 L_1 r_2 \sin q_2 \dot{q}_1$, $B_{21} = -2m_2 L_1 r_2 \sin q_2 \dot{q}_1$, $B_{22} = 0$, $G_2 = m_2 r_2 g \cos(q_1 + q_2)$, $G_1 = (m_1 r_1 + m_2 L_1) g \cos q_1 + m_2 r_2 g \cos(q_1 + q_2)$.

The nominal arguments for the system are selected as $m_1 = m_2 = 1kg$, $L_1 = L_2 = 0.2m$, $r_1 = r_2 = 0.1m$, $I_1 = 0.64kg \cdot m^2, I_2 = 0.16kg \cdot m^2$ and $g = 9.8m/s^2$.

Compensation of contour tracking control and trajectory tracking control are carried out in task space and joint space, respectively, so it is necessary to transform the control space between task space and joint space. The relationship between task space and joint space error is established as follows:

$$e_c = J(q)e_q \tag{44}$$

where $e_c$ is tracking error vector of the end-effector in task space, $e_q$ denotes the tracking error vector in joint space. $J(q)$ is Jacobian matrix of two-link robot manipulator which is defined as follows

$$J(q) = \begin{bmatrix} -L_1 sin(q_1) - L_2 sin(q_1 + q_2) & -L_2 sin(q_1 + q_2) \\ L_1 cos(q_1) + L_2 cos(q_1 + q_2) & L_2 cos(q_1 + q_2) \end{bmatrix}$$

By using (40), the contour error can be derived as

$$\varepsilon = \begin{bmatrix} c_x & c_y \end{bmatrix} e_c = \begin{bmatrix} c_x & c_y \end{bmatrix} J(q) e_q \tag{45}$$

The contouring error compensation rectifier gain of two-link robotic manipulator can be derived as

$$\begin{bmatrix} C_1 & C_2 \end{bmatrix} = \begin{bmatrix} c_x & c_y \end{bmatrix} J(q) \tag{46}$$

Considering the influence of the actuator faults and the lumped uncertainty of the system on control performances of robotic manipulator, $\varphi(q,\dot{q},t)$ and $\phi(q,\dot{q},\tau)$ are selected as follows

$$\varphi(q,\dot{q},t) = \begin{bmatrix} 0.5\dot{q}_1 + sin(3q_1) + 0.5 sin(\dot{q}_1) \\ 1.3\dot{q}_2 - 1.8 sin(2q_2) + 1.1 sin(\dot{q}_2) \end{bmatrix}$$

and $\phi(q,\dot{q},\tau) = \begin{cases} 30 sin(q_1 q_2) + 4 cos(\dot{q}_1 q_2) + 15 cos(\dot{q}_1 \dot{q}_2) & T_{f1} \geq 1.5 \\ 0 & T_{f2} \geq 1.5 \end{cases}$

A. Simulation Results of Path Tracking without Planning

In the situation of fast-speed motion of the robotic manipulator, the simulation experiment is carried out. The simulation time is set to $3s$, sampling step is $0.001s$, the premier state of system is $[0.3 \quad 0.3]^T$, the x-axis and y-axis trajectories are designed as follows

$$x = \begin{cases} 3 - 2t, & 0 \leq t \leq 1 \\ 1, & 1 < t \leq 2 \\ 2(t-1) - 1, & 2 < t \leq 3 \end{cases}, \quad y = \begin{cases} 2t + 1, & 0 \leq t \leq 1 \\ 3 - 2(t-1), & 1 < t \leq 2 \\ 1, & 2 < t \leq 3 \end{cases} \tag{47}$$

The traditional PID controller parameters are given by $K_P = diag[4500,4500]$, $K_I = diag[150,150]$, $K_D = diag[650,650]$. The parameter values used in the ANTSMC are $c_1 = diag[10,10]$, $c_2 = diag[10,13.5]$, $\eta = diag[0.5,0.5]$, $\alpha = 3$, $\beta = 5$, $\upsilon = diag[50,165]$, $\xi = 10$, and the proportional gain and differential gain of the cross coupling controller (CCC) are $K_p = 500$, $K_d = 100$, respectively.

Trajectory tracking response of joint 1 and 2 in joint space are given in Fig 5. The error of trajectory tracking of joint 1 and 2 in joint space are shown in Fig.6. Corresponding to joint space, the contour tracking of the end-effector of robotic manipulator is displayed in Fig.7.

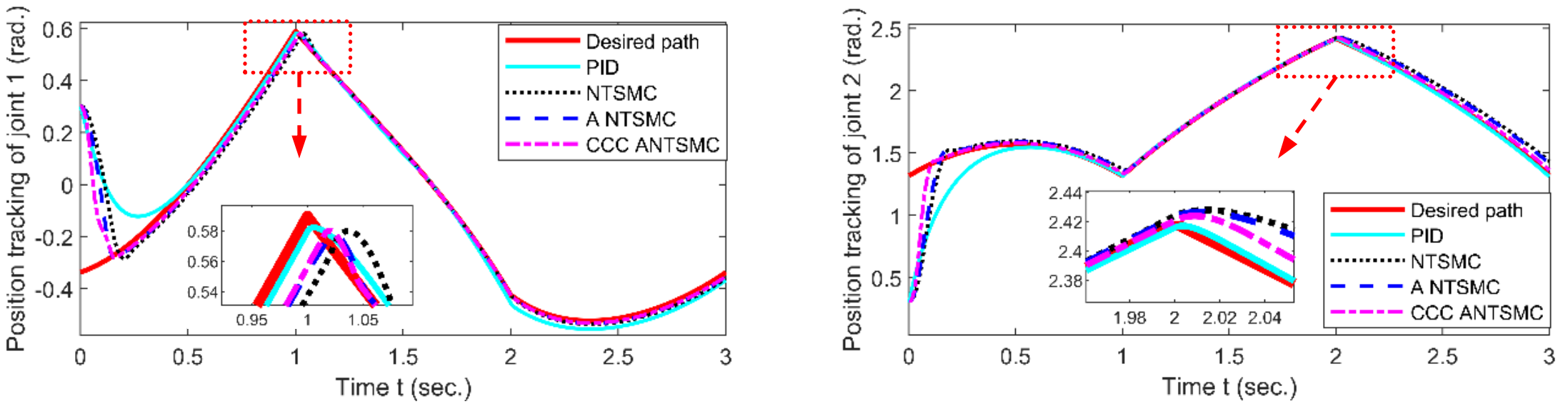

(a) joint 1 (b) joint 2

Fig.5 Trajectory tracking of joint 1 and 2 in joint space

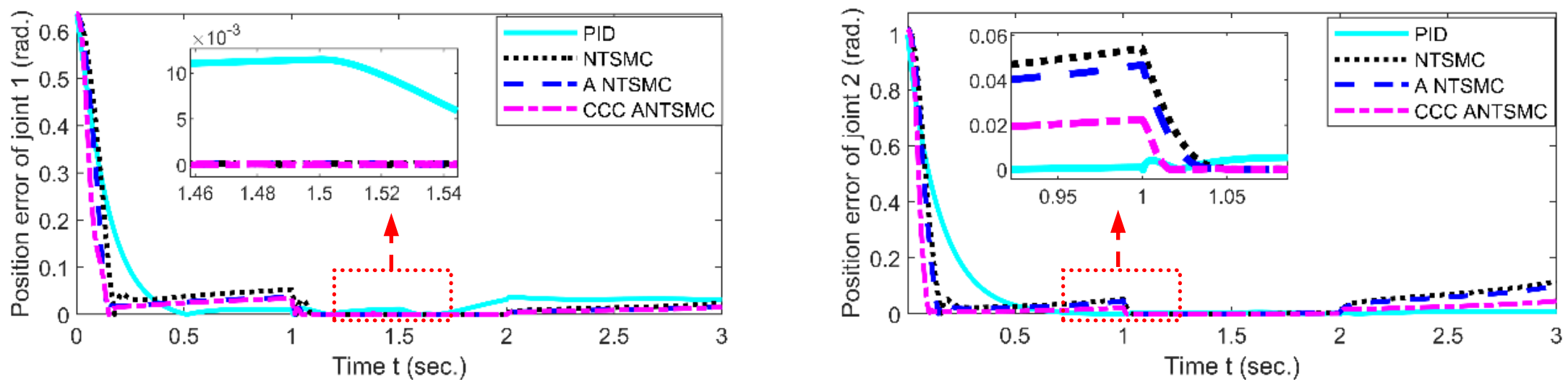

(a) joint 1 (b) joint 2

Fig.6 Trajectory tracking error of joint 1 and 2 in joint space

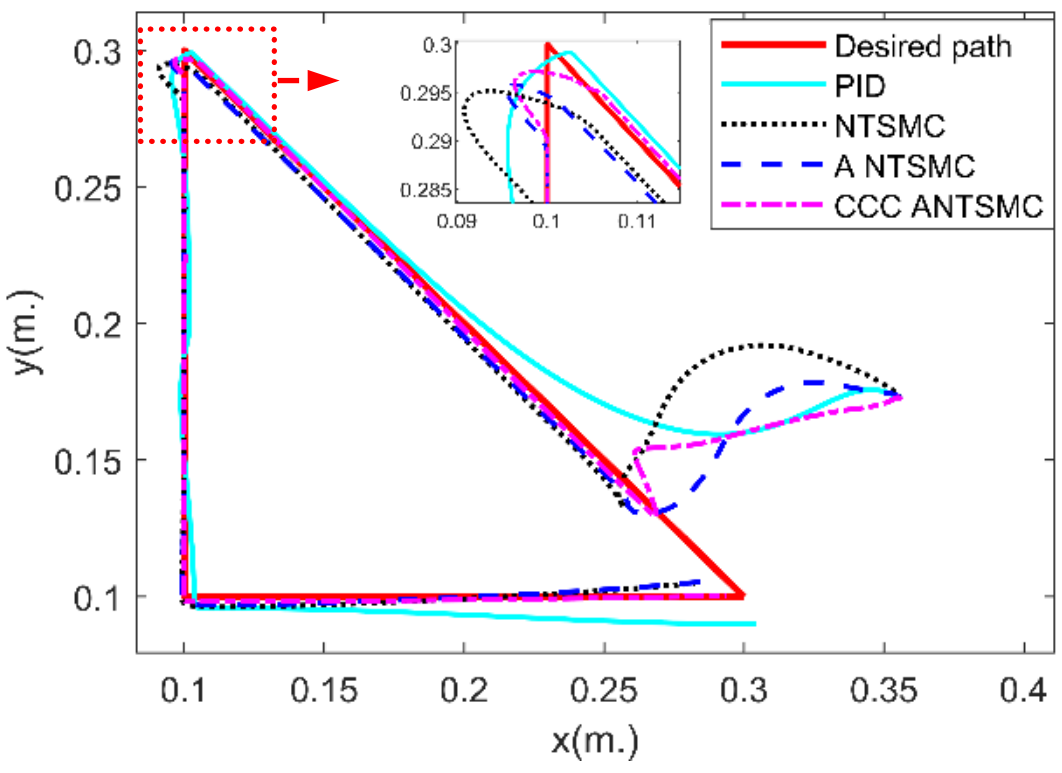

Fig.7 Contour tracking of end-effector of robotic manipulator

For a clear comparison, the average tracking error can be defined as

$$e_j = \sqrt{\frac{1}{N}\sum_{i=1}^{N}(\|e_j(k)\|^2)} \qquad j = 1,2$$

where $N$ is the number of simulation steps.

The performance comparison of each controller is given by Table 1.

Table 1 Performance comparison of each controller

| Controller | Average tracking error of joint 1/(rad) | Average tracking error of joint 2/(rad) | Average error of contour/(m) |
|---|---|---|---|
| PID | 0.0461 | 0.0826 | 0.0076 |
| NTSMC | 0.0467 | 0.0715 | 0.0051 |
| A NTSMC | 0.0387 | 0.0651 | 0.0044 |
| CCC NTSMC | 0.0297 | 0.0482 | 0.0029 |

Simulation results demonstrate that the proposed control strategy can optimize the contour error in the initial stage quickly, and the tracking performance of each joint is also improved. However, due to the sudden change of speed and acceleration at the inflection point of the adjacent line segment, the control performance of each controller will be weakened to a certain extent.

B. Simulation Results of Path Tracking via Parabola Transition Planning

To analyze the influence of position and speed mutation at inflection point of non-smooth path in workspace on control performance of robotic manipulator, the parameters of each controller are consistent with those in the previous section. Parabolic transition planning technology is used to smooth transition between adjacent linear segments, and the interpolation points of linear segments are set as follows

$$x = x_0 + \Pi \Delta x, \quad y = y_0 + \Pi \Delta y$$

where $x_0$, $y_0$ are the starting position of the line, $\Delta x$ and $\Delta y$ are position increments along with $x$ and $y$ axis, $\Pi$ is normalization factor [27].

Trajectory tracking response of joint 1 and 2 in joint space are given in Fig 8. The error of trajectory tracking of joint 1 and 2 in joint space are shown in Fig 9. Corresponding to joint space, the contour tracking of the end-effector of robotic manipulator is given in Fig 10.

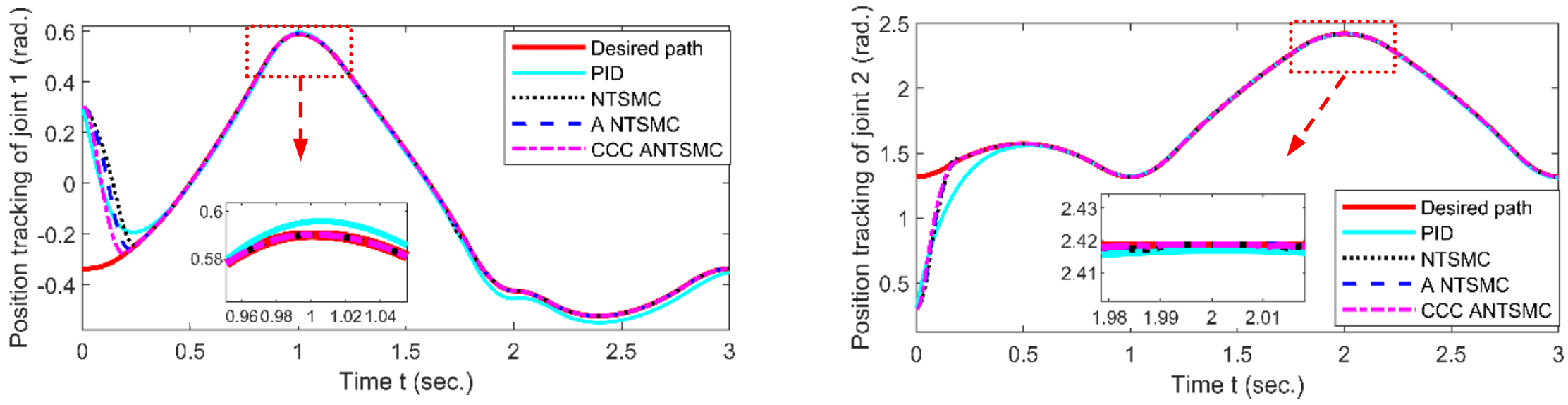

(a) joint 1 (b) joint 2

Fig.8 Trajectory tracking of joint 1 and 2 in joint space (Parabola Transition Planning)

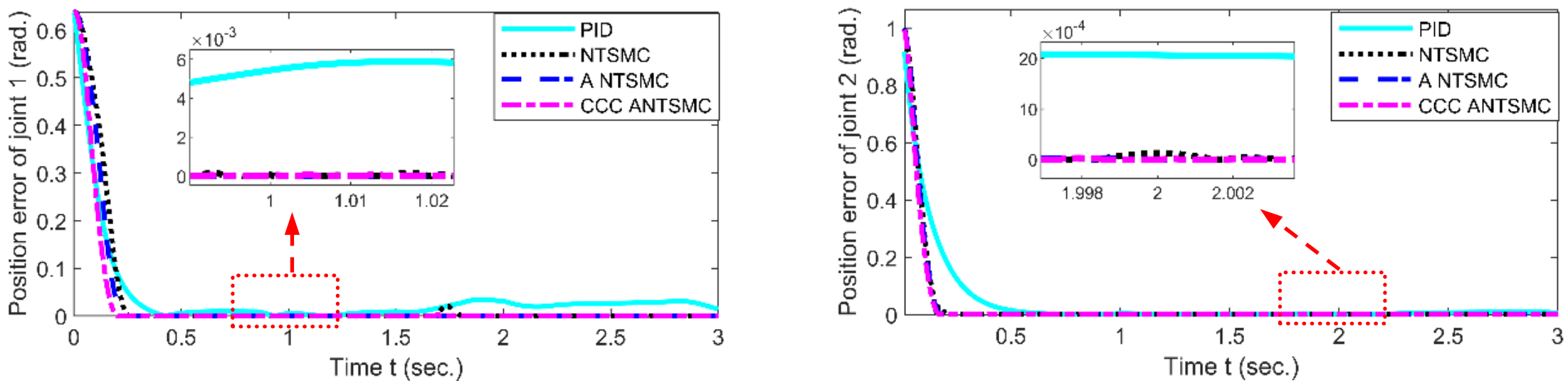

(a) joint 1 (b) joint 2

Fig.9 Trajectory tracking error of joint 1 and 2 in joint space (Parabola Transition Planning)

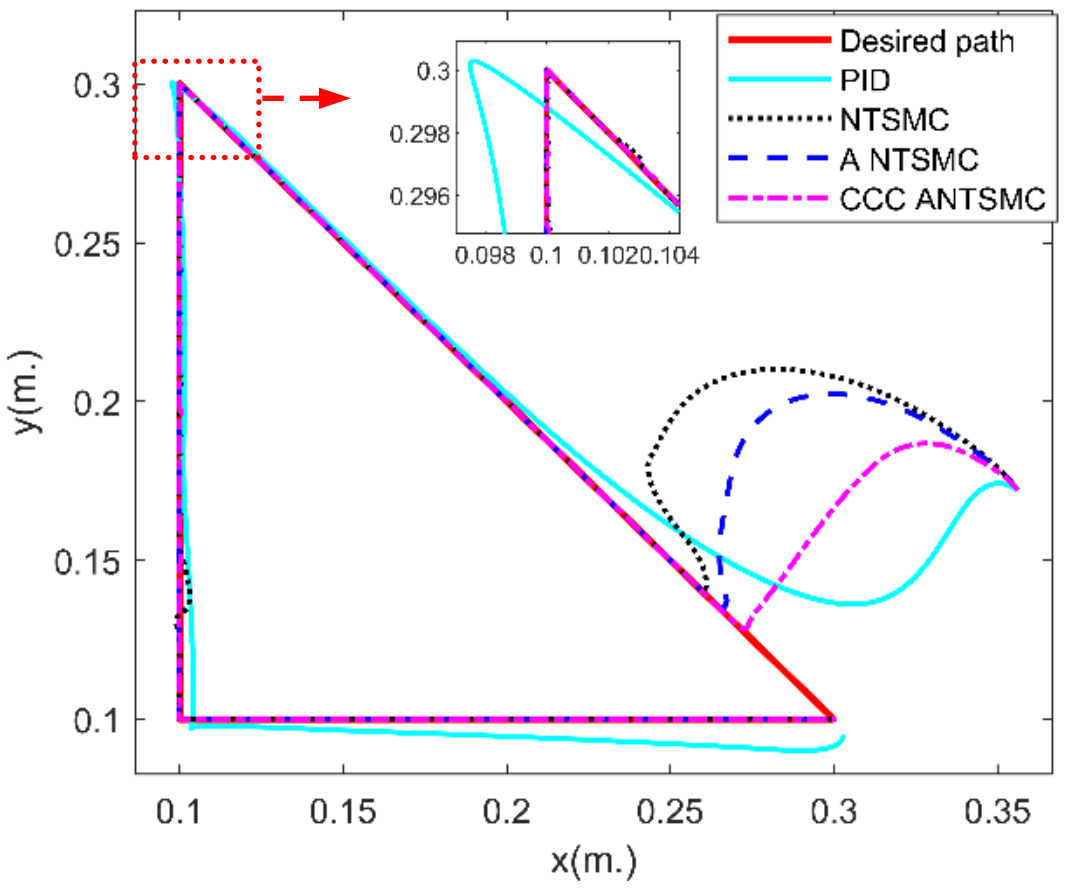

Fig.10 Contour tracking of end-effector of robotic manipulator (Parabola Transition Planning)

The performance comparison of each controller with parabola transition planning is given by Table 2.

Table 2 Performance comparison of each controller (Parabola Transition Planning)

| Controller | Average tracking error of joint 1/(rad) | Average tracking error of joint 2/(rad) | Average error of contour/(m) |
|---|---|---|---|
| PID | 0.0353 | 0.0585 | 0.0045 |
| NTSMC | 0.0276 | 0.0289 | 0.0024 |
| A NTSMC | 0.0223 | 0.0269 | 0.0023 |
| CCC NTSMC | 0.0145 | 0.0138 | 0.0011 |

The simulation results show that the control property of each control strategy will be optimized to a certain extent at the connection inflection point of the straight line segment through the path planning of parabolic transition. Furthermore, Fig.9 and Fig. 10 show that the designed ANTSMC based on cross coupling has better contour tracking performance when the dynamic performance of two joints is different. Compared without the cross coupled control method, the contour tracking error of end-effector of the designed control system is smaller, and it has better robustness to dynamic uncertainties, actuator faults and external disturbances.

## 5.2. Experiments

In this part, the experimental platform includes three parts as shown in Fig. 11: the upper system, the lower system and experimental prototype of robotic manipulator. Among them, the upper system is a hardware and software system composed of PC and C++/MFC to communicate with the lower system. The lower computer adopts the four-axis motion control card (GTS-400-PV-PCI) of Googoltech Company, and forms a closed-loop control system with the data acquisition module such as encoder and analog output module and the experimental prototype.

The experimental platform is designed as a 2-DOF robotic manipulator, the length of limb 1 is $0.2m$, the length of limb 2 is $0.2m$, the two joints are driven by two sets of AC servo motors, the joint 1 is driven by a Delta servo motor, and the joint 2 is driven by a Panasonic servo motor. Both joints adopt belt transmission device to transfer torque, and joint angular displacement and angular velocity information are detected by an encoder with a resolution of 2500p/r.

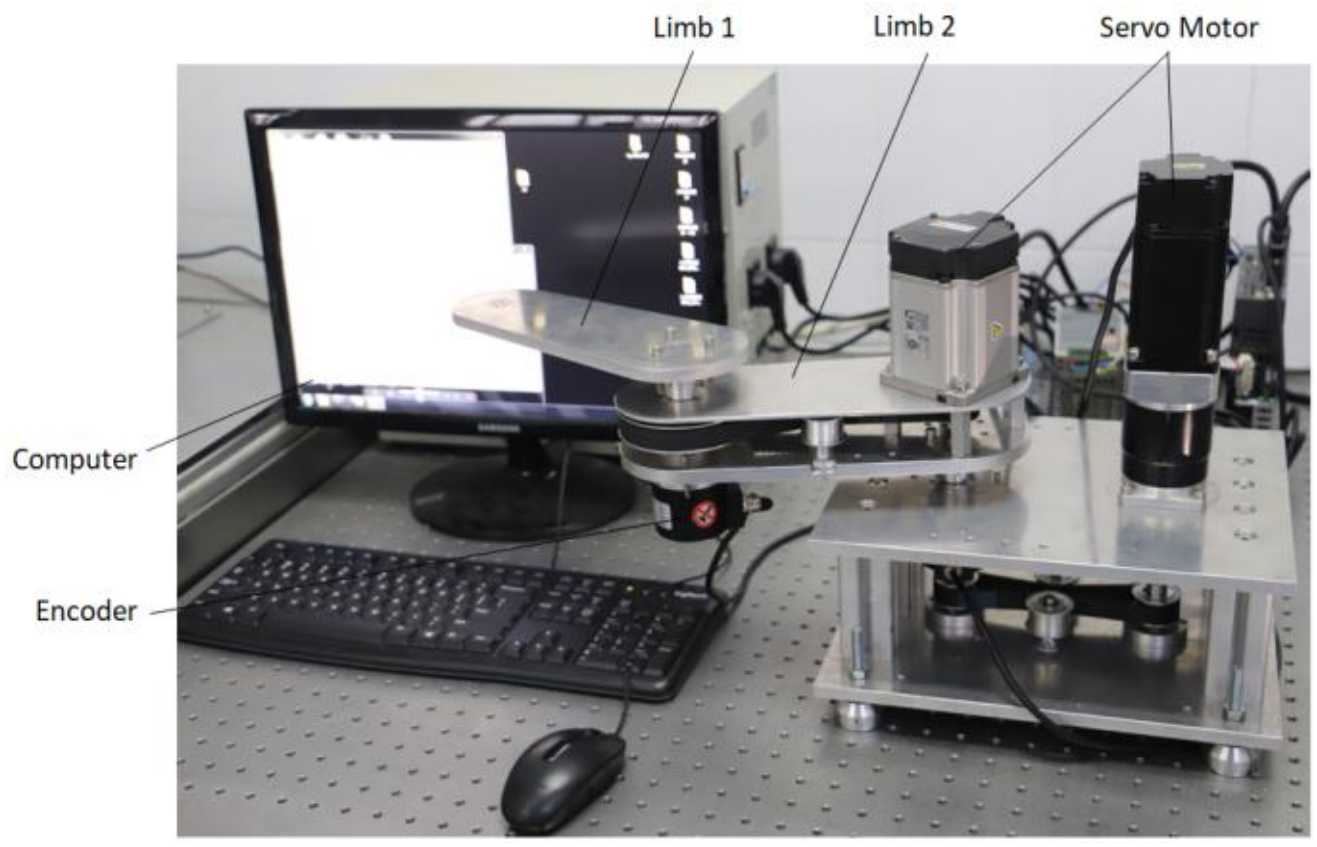

Fig.11 The experimental platform of robotic manipulator

The experiment time is set as 3s, and the number of sampling points is 30001. The traditional PID controller parameters are given by $K_P = diag[3.5,0.8]$, $K_I = diag[0.09,0.012]$, $K_D = diag[0.15,0.1]$. The parameter values used in the ANTSMC are $c_1 = diag[12.05,12.02]$ , $c_2 = diag[0.5,5.5]$ , $\eta = diag[0.5,0.5]$, $\alpha = 3$, $\beta = 5$, $\upsilon = diag[40.5,42.2]$, $\xi = 5$, and the proportional gain and differential gain of the cross coupling controller (CCC) are $K_p = 0.5$, $K_d = 0.02$, respectively.

Trajectory tracking response of joint 1 and 2 in joint space are given in Fig 12. Corresponding to joint space, the contour tracking of the end-effector of robotic manipulator is given in Fig.13.

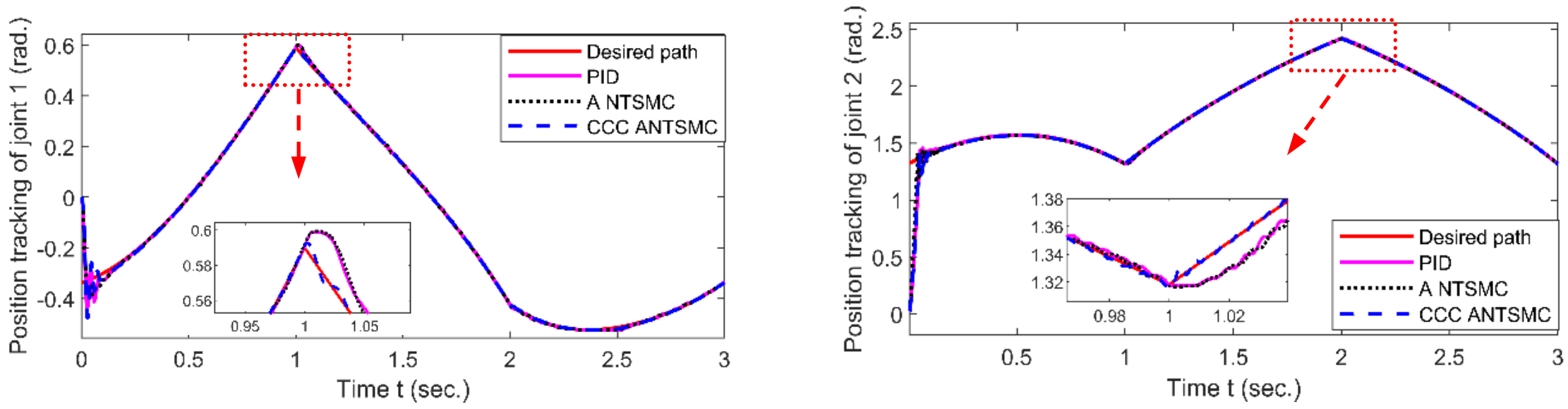

(a) joint 1 (b) joint 2

Fig.12 Trajectory tracking of joint 1 and 2 in joint space

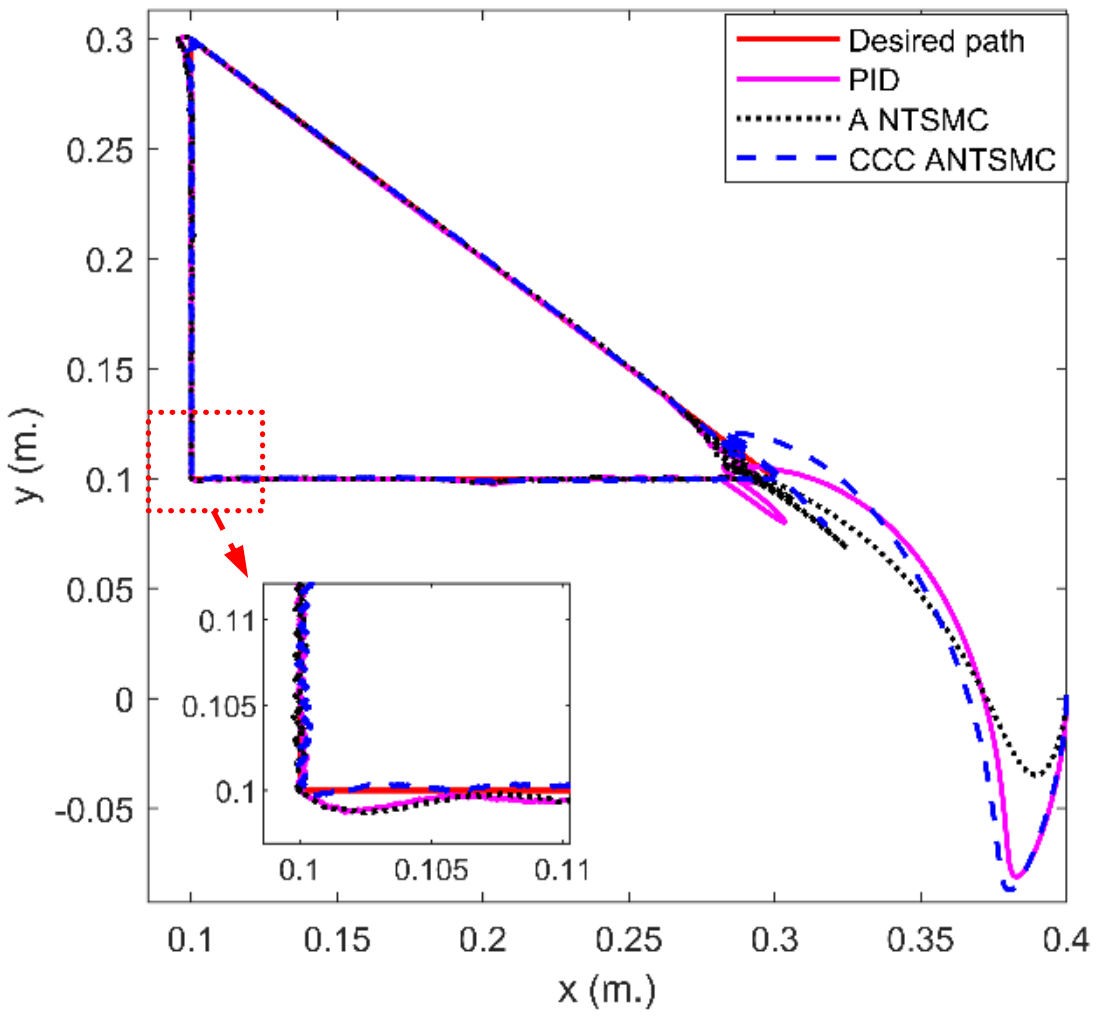

Fig.13 Contour tracking of end-effector of robotic manipulator

The performance comparison of each controller is given by Table 3.

Table 3 Performance comparison of each controller

| Controller | Average tracking error of joint 1/(rad) | Average tracking error of joint 2/(rad) | Average error of contour/(m) |
|---|---|---|---|
| PID | 0.0415 | 0.0439 | 0.0062 |
| A NTSMC | 0.0304 | 0.0431 | 0.0053 |
| CCC NTSMC | 0.0250 | 0.0269 | 0.0028 |

Experimental results demonstrate that the proposed control strategy has better property than the PID control method in profile control. Due to sudden changes in position, velocity and acceleration at the junction inflection point of adjacent straight segments, each controller will suffer performance degradation to a certain extent at the time of 1s and 2s, but the adaptive tuning method of the proposed control strategy for this kind of mutation is better than the traditional PID controller and the contour tracking error of the cross coupling compensation control strategy is small.

To analyze the influence of position and speed mutation at inflection point of non-smooth path in workspace on robot performance, the parameters of the above controller keeps consistent. Trajectory tracking response of joint 1 and 2 in joint space are given in Fig 14. Corresponding to joint space, the contour tracking of the end-effector of robotic manipulator is given in Fig.15.

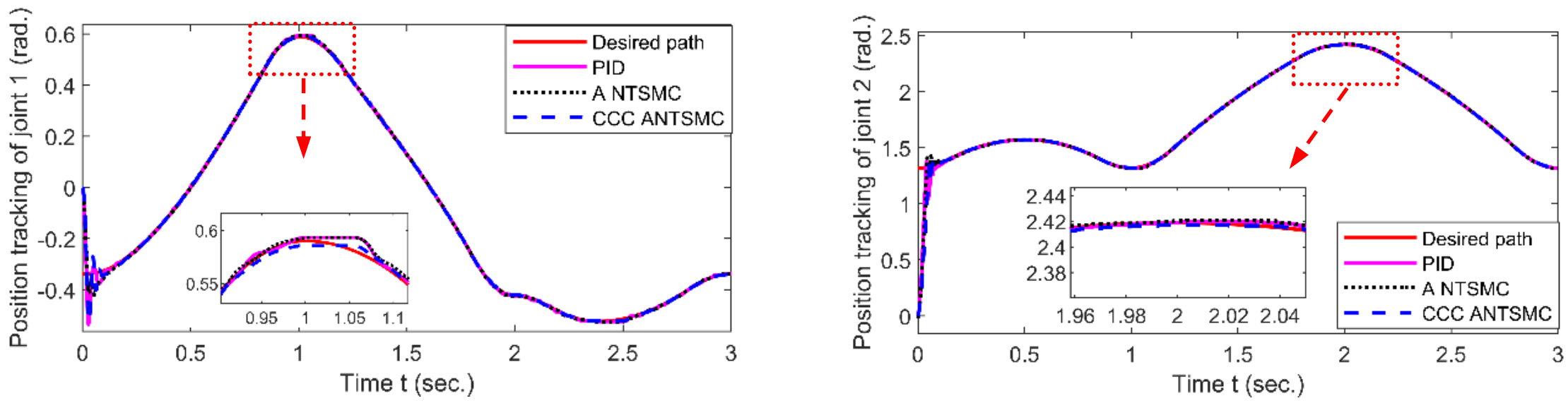

(b) joint 1 (b) joint 2

Fig.14 Trajectory tracking of joint 1 and 2 in joint space (Parabola Transition Planning)

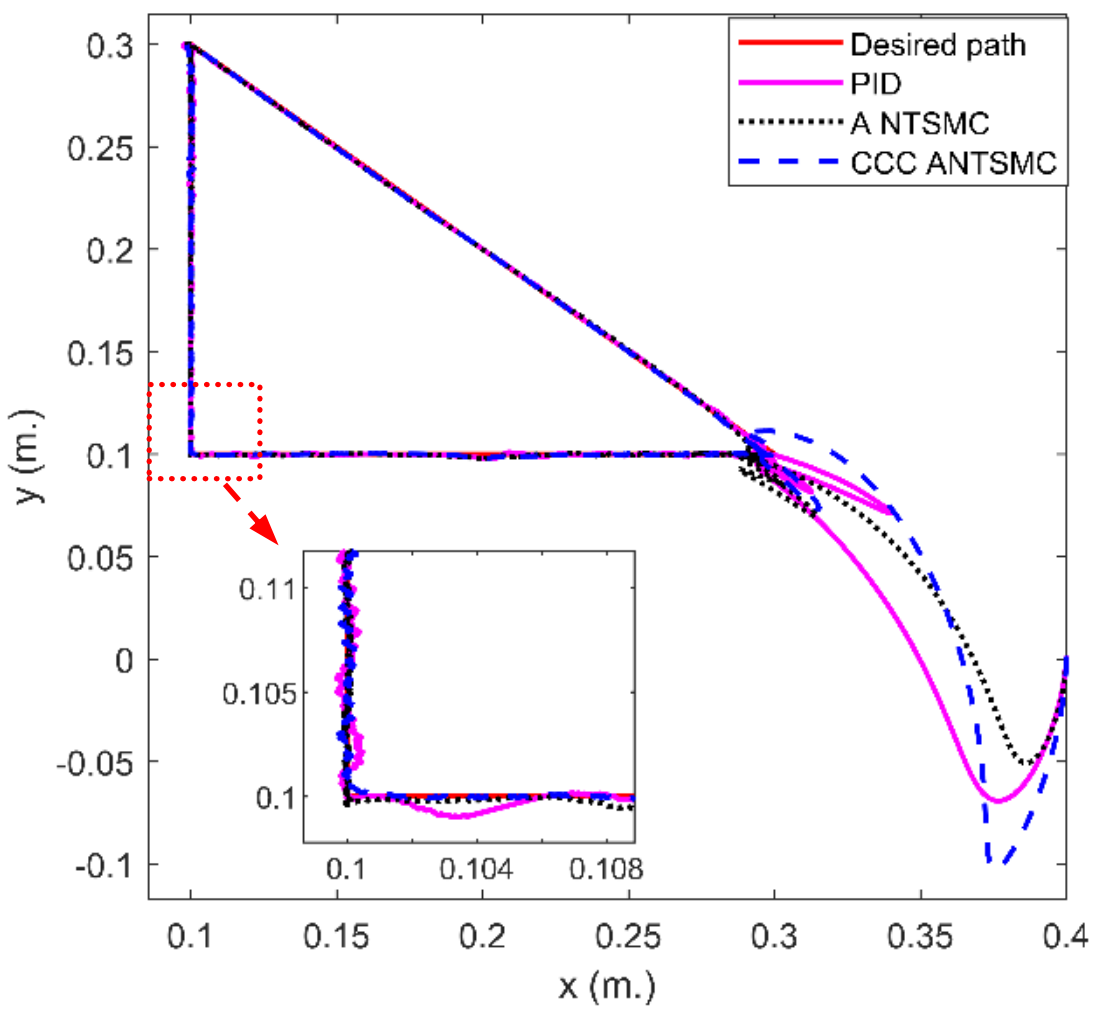

Fig.15 Trajectory tracking of end-effector of robotic manipulator (Parabola Transition Planning)

The performance comparison of each controller with parabola transition planning is shown by Table 4

Table 4 Performance comparison of each controller (Parabola Transition Planning)

| Controller | Average tracking error of joint 1/(rad) | Average tracking error of joint 2/(rad) | Average error of contour/(m) |
|---|---|---|---|
| PID | 0.0232 | 0.0452 | 0.0046 |
| A NTSMC | 0.0215 | 0.0335 | 0.0041 |
| CCC NTSMC | 0.0224 | 0.0101 | 0.0019 |

Experimental results reveal that the proposed control strategy has better control property than the traditional PID control method. After the parabolic transition path planning, the control performance of each controller will be improved to a certain extent at inflection point of non-smooth path in workspace. It can be seen from Fig. 14 and Fig.15 that the end-effector of robotic manipulator has better contour tracking property when the dynamic performance of the two joints is different. Compared with the decoupled control method, the proposed controller has smaller contour tracking error and advantageous robustness to dynamic uncertainties.

It is worth noting that the belt transmission device is used to transfer torque in both joints of the experimental platform, so the belt tensioning process is required at the start time, which affects the control performance of robotic manipulator to a certain extent in the initial stage.

## 6. Conclusions

In this paper, an ANTSMC law based on cross coupling is proposed for the issue of precis contour tracking of robotic manipulator under dynamic uncertainties, coupling and actuator faults. This adaptive tuning approach is selected based on the strengths of the NTSMC and can drive the system state to the equilibrium point in finite time. While the proposed ANTSMC can improve the trajectory tracking accuracy of each joint of the robotic manipulator, coupling factors are introduced among the multi-axes to optimize the coordination control property of robotic manipulator. A unified framework of cross-coupling contour compensation and reference position pre-compensation is designed by combining cross-coupling control technology with parabolic transition trajectory planning technology. Simulation and experimental results illustrate that the proposed control strategy has better tracking property and by improving the coordination of each joint, the influence of system uncertainty on the contour tracking accuracy of the end-effector of robot is reduced, and the accurate contour tracking control of the end-effector on the preset machining path is realized, which can achieve higher contour control accuracy compared with the traditional trajectory control method.

## Funding

This work was supported by the National Natural Science Foundation of China (Grant numbers. 51905115), and Scientific research project of Guangzhou Education Bureau, China (Grant numbers. 202032821).

## Conflicts of Interest/Competing interests

The authors declared that they have no conflicts of interest to this work.

## Authors' Contributions

All authors contributed to the study conception and design. Material preparation, data collection and analysis were performed by Du Baolin, Cui Aodong and Zhu Puchen. The first draft of the manuscript was written by Zhu Dachang and all authors commented on previous versions of the manuscript. All authors read and approved the final manuscript.